\documentclass[12pt]{amsart}

\usepackage{amsfonts,amssymb,amsthm,amsmath}

\newtheorem{theorem}{Theorem}[section]
\newtheorem{lemma}[theorem]{Lemma}
\newtheorem{proposition}[theorem]{Proposition}
\newtheorem{corollary}[theorem]{Corollary}

\theoremstyle{definition}
\newtheorem{definition}[theorem]{Definition}

\newtheorem{conjecture}[theorem]{Conjecture}
\newtheorem*{acknowledgments}{Acknowledgments}

\newcommand{\C}{{\mathbb{C}}}
\newcommand{\Q}{{\mathbb{Q}}}

\newcommand{\N}{{\mathbb{N}}}
\renewcommand{\P}{{\mathbb{P}}}
                                       
\begin{document}

\title{Effective freeness of adjoint line bundles}

\author[Gordon Heier]{Gordon Heier}

\address{Ruhr-Universit\"at Bochum\\
Fakult\"at f\"ur Mathematik\\
D-44780 Bochum\\
Germany}

\curraddr{Harvard University\\
Department of Mathematics\\
Cambridge, MA 02138\\
USA}

\email{heier@cplx.ruhr-uni-bochum.de}

\date{August 17, 2001.}



\begin{abstract}
This note shows how two existing methods to provide effective bounds for the freeness of adjoint line bundles can be applied in combination to establish a new such bound which (approximately) differs from the linear bound conjectured by Fujita only by a factor of the cube root of the dimension of the underlying manifold. As an application, a new effective statement for pluricanonical embeddings is derived.
\end{abstract}

\maketitle

\tableofcontents

\section{Introduction and Statement of the Main Theorem}
Let $L$ be an ample line bundle over a compact complex projective manifold $X$ of complex dimension $n$. Let $K_X$ be the canonical line bundle of $X$. The following conjecture is due to Fujita \cite{F}.
\begin{conjecture}[Fujita]
The adjoint line bundle $K_X+mL$ is base point free (i.e. spanned by global holomorphic sections) for $m\geq n+1$. It is very ample for $m\geq n+2$.
\end{conjecture}
The standard example of $X=\P^n$ shows that the conjectured numerical bounds are optimal in the sense that there is a counter example to the statement of the conjecture if the bounds are lowered. In the case of $X$ being a compact Riemann surface, the conjecture is easily verified by means of the Riemann-Roch theorem. Moreover, Reider \cite{R} was able to validate the conjecture also in the case $n=2$. In higher dimensions, the very ampleness part of the conjecture has proved to be quite intractable so far. In fact, no further results seem to be known here. On the other hand, several further results have been established towards the freeness conjecture. The case $n=3$ was solved by Ein and Lazarsfeld \cite{EL} (see also \cite{Fujpre}), and $n=4$ is due to Kawamata \cite{K}. In arbitrary dimension $n$, the state-of-the-art is that $K_X+mL$ is base point free for any integer $m$ which is no less than a number roughly of order $n^2$ (see below for exact statements).\par
To the author's knowledge, \cite{DBourbaki} constitutes the most recent survey on the subject under discussion. It contains an extensive list of references (see also the references at the end of this article) and, furthermore, introduces the reader to various other effective results in algebraic geometry.\par
The above-mentioned bound in the case of arbitrary dimension $n$ can be derived from each of the following two theorems due to Angehrn and Siu \cite{AS} (see also \cite{S}) and Helmke \cite{H}, \cite{H2}, respectively. Their proofs essentially adhere to the same inductive strategy, and we will combine the respective advantages of these two methods to obtain a value for the freeness bound which is even closer to the conjectured linear bound.\par
First, let us state the bound given by \cite{AS}.
\begin{theorem} [\cite{AS}]
The line bundle $K_X+mL$ is base point free for $m \geq \frac{1}{2}n(n+1)+1$.
\end{theorem}
Secondly, we state Helmke's result. Due to the nature of his technique, the assumptions of his theorem are formulated in a slightly different way. We quote the result in the way it is presented in \cite{H}, because the slight improvement achieved in \cite{H2} is not relevant for our purposes.
\begin{theorem} [\cite{H}] \label{helmkethm}
Assume that $L$ has the additional properties that
\begin{equation*} L^n> n^n\end{equation*}
and for all $x\in X$:
\begin{equation*}L^d.Z\geq m_x(Z)\cdot n^d\end{equation*}
for all subvarieties $Z \subset X$ with $x\in Z,\ d= \dim Z \leq n-1$ and multiplicity $m_x(Z)\leq{{n-1}\choose{d-1}}$ at $x$. Then $K_X+L$ is base point free.
\end{theorem}
If $n\geq 3$, it is clear that we need to set
\begin{equation*}m_0:= \max\{n\cdot\sqrt[d]{{n-1\choose{d-1}}}:d\in \N \text{ and }1\le d\le n\}\end{equation*}
in order to determine the minimal bound $m_0$ deducible from Theorem \ref{helmkethm} such that $K_X+mL$ is base point free for any integer $m \geq m_0$. Since 
\begin{equation*}\sqrt[d]{{n-1\choose{d-1}}} \geq \frac{1}{\sqrt[d]3}\frac{n^{1-\frac{1}{d}}}{d}\end{equation*}
according to our Lemma \ref{coeffestim}, we find that the $m_0$ which one can derive from Theorem \ref{helmkethm} is essentially also of the order $n^2$.\par
We conclude this section with the statement of our Main Theorem, which asserts that the bound $m_0$ can be chosen to be a number of the order $n^\frac 4 3$. For its proof, see Section \ref{sectionmainthmproof}.
\begin{theorem}[Main Theorem]
The line bundle $K_X+mL$ is base point free for any integer $m$ with
\begin{equation*}
m\geq(e+\frac 1 2)n^\frac 4 3 + \frac 1 2 n^\frac 2 3 +1,
\end{equation*}
where $e\approx 2.718$ is Euler's number.
\end{theorem}

\section{Estimates for Binomial Coefficients}
In order to precisely understand the nature of the numerical conditions in the assumptions of Theorem \ref{helmkethm}, we prove some auxiliary estimates in this section. We begin with the following lemma.
\begin{lemma}\label{l'hopital}
For all $x \in \,]0,1[\,:\ 1<\left(\frac{1}{1-x}\right)^{\frac{1-x}{x}}<e.$
\end{lemma}
\begin{proof}
It is obvious that $1$ is a strict lower bound of the given expression, so it remains to show that
\begin{equation*}\left(\frac{1}{1-x}\right) ^{\frac{1-x}{x}}<e.\end{equation*}
Taking $\log$ on both sides of the inequality, we see that we are done if we can show that
\begin{equation*}g(x):=\frac{x-1}{x}\log(1-x)<1\end{equation*}
on the open unit interval. However, for this it suffices to prove that $\lim_{x\rightarrow 0^+}g(x)=1$ and $g'(x)<0$. The former is easily verified using L'H\^ opital's rule, while the latter follows readily from a simple computation.
\end{proof}
In the proof of the subsequent Lemma \ref{coeffestim}, we will employ Lemma \ref{l'hopital} in the form of the following corollary.
\begin{corollary}\label{estim}
Let $n$ be an integer $\geq 2$. Let $d$ be an integer with $1 \leq d \leq n-1$. Then
\begin{equation*}1<\left(\frac{n}{n-d}\right)^\frac{n-d}{d}<e.\end{equation*}
\end{corollary}
\begin{proof}
We have
\begin{eqnarray*}
\left(\frac{n}{n-d}\right)^\frac{n-d}{d}=\left(\frac{1}{1-\frac{d}{n}}\right)^\frac{1-\frac{d}{n}}{\frac{d}{n}}.
\end{eqnarray*}
Thus the corollary follows immediately from Lemma \ref{l'hopital}.
\end{proof}
The preceding considerations allow us to estimate the binomial coefficients from Theorem \ref{helmkethm} in the form of the following lemma, which is the key ingredient in our numerical considerations.
\begin{lemma}\label{coeffestim}
Let $1\leq d \leq n-1$. Then
\begin{equation*}\frac{1}{\sqrt[d]3}\frac{n^{1-\frac{1}{d}}}{d}\leq \sqrt[d]{{n-1\choose{d-1}}} \leq e\frac{n}{d}.\end{equation*}
\end{lemma}
\begin{proof}
In \cite{A}, page 206, Stirling's formula is stated as
\begin{equation*}
\Gamma(x)=\sqrt{2\pi}\,x^{x-{\frac 1 2}}e^{-x}e^{\frac {\theta(x)} {12x}}
\end{equation*}
for $x>0$ with $0<\theta(x)<1$. In particular,
\begin{equation*}
\sqrt{2\pi}\,x^{x-{\frac 1 2}}e^{-x}\leq\Gamma(x)\leq\sqrt{2\pi}\,x^{x-{\frac 1 2}}e^{-x}e^{\frac 1 {12}}
\end{equation*}
for any $x\geq 1$. Thus, for the proof of the desired estimate from above, Stirling's formula enables us to proceed as follows.
\begin{eqnarray*}
&& {n-1\choose d-1}={\frac{(n-1)!}{(d-1)!(n-d)!}} ={\frac 1
{n-d}}{\frac{\Gamma(n)}{\Gamma(d)\Gamma(n-d)}}\\
&\leq& {\frac 1 {n-d}}\frac {\sqrt{2\pi}\,n^{n-{\frac 1 2}}e^{-n}e^{\frac 1 {12}}}
{\sqrt{2\pi}\,d^{d-{\frac 1 2}}e^{-d}
\sqrt{2\pi}\,(n-d)^{n-d-{\frac 1 2}}e^{-(n-d)}}\\
&=&\frac{e^{\frac 1 {12}}}{\sqrt{2\pi}}\sqrt{{\frac d{(n-d)n}}} \left({\frac n d}\right)^d\left(\frac n {n-d}\right)^{n-d} \\
&\leq& \left({\frac n d}\right)^d\left({\frac n {n-d}}\right)^{n-d}.
\end{eqnarray*}
Resorting to Corollary \ref{estim}, we eventually conclude:
\begin{equation*}
{n-1\choose d-1}^{\frac 1 d}\leq{\frac n d}\left({\frac n {n-d}}\right)^{\frac{n-d} d}\leq e \frac n d.
\end{equation*}\par
The desired estimate from below is proved analogously:
\begin{eqnarray*}
&& {n-1\choose d-1}={\frac{(n-1)!}{(d-1)!(n-d)!}} ={\frac 1{n-d}}{\frac{\Gamma(n)}{\Gamma(d)\Gamma(n-d)}}\\
&\geq&{\frac 1{n-d}}{\frac{\sqrt{2\pi}\,n^{n-{\frac 1 2}}e^{-n}}{\sqrt{2\pi}\,d^{d-{\frac 1 2}}e^{-d}e^{\frac 1 {12}} \sqrt{2\pi}\,(n-d)^{n-d-{\frac 1 2}}e^{-(n-d)}e^{\frac 1 {12}}}}\\
&=&\frac 1 {e^{\frac 1 6}\sqrt{2\pi}}\sqrt{{\frac d{(n-d)n}}} \left({\frac n d}\right)^d\left({\frac n{n-d}}\right)^{n-d} \\
&\geq&\frac 1 3 \frac 1 n  \left({\frac n d}\right)^d\left({\frac n {n-d}}\right)^{n-d} \\
\end{eqnarray*}
Using Corollary \ref{estim} again, we obtain:
\begin{equation*}
{n-1\choose d-1}^{\frac 1 d}\geq\sqrt[d]{\frac 1 {3n}}\frac n d=\frac{1}{\sqrt[d]3}\frac{n^{1-\frac{1}{d}}}{d}.
\end{equation*}
\end{proof}

\section{Proof of the Main Theorem}\label{sectionmainthmproof}
The following theorem states the improved effective freeness bound which we shall prove at the end of this section.
\begin{theorem}[Main Theorem]\label{mainthm}
The line bundle $K_X+mL$ is base point free for any integer $m$ with
\begin{equation*}
m\geq(e+\frac 1 2)n^\frac 4 3 + \frac 1 2 n^\frac 2 3 +1,
\end{equation*}
where $e\approx 2.718$ is Euler's number.
\end{theorem}
First of all, let us recall how a result of this type can be proved by means of multiplier ideal sheaves.\par
 Let $x\in X$ be an arbitrary but fixed point. The key idea of both \cite{AS} and \cite{H} is to find an integer $m_0$ (as small as possible) and a singular metric $h$ of the line bundle $m_0L$ with the following two properties:
\begin{enumerate}
\item \label{prop1}Let $h$ be given locally by $e^{-\varphi}$. Then the curvature current $i\partial\bar\partial\varphi$ dominates a positive definite smooth $(1,1)$-form on $X$ in the sense of currents.
\item Let the multiplier ideal sheaf of $h$ be defined stalk wise by 
\begin{equation*}(\mathcal I_h)_x:=\{f\in \mathcal O_{X,x}:|f|^2e^{-\varphi} \text{ is locally integrable at }x\}.\end{equation*}
Then, in a neighborhood of $x$, the zero set of $\mathcal I_h$, which we denote by $V(\mathcal I_h)$, is just the point $x$. (This is the key property we are looking for. Note that the support of $V(\mathcal I_h)$ is just the set of points where $h$ is not locally integrable.)
\end{enumerate}
The first property implies that 
\begin{equation*}H^q(X,\mathcal I_h(K_X+m_0L))=(0) \quad (q\geq 1),\end{equation*}
due to the vanishing theorem of Nadel \cite{N}, \cite{N2}. (In the special case when the singular metric is algebraic geometrically defined, Nadel's vanishing theorem is the same as the theorem of Kawamata and Viehweg \cite{K2}, \cite{V}.)
With this information and the second property, it is easy to obtain an element of $\Gamma(X,K_X+m_0L)$ which does not vanish at $x$. Namely, consider the the standard short exact sequence
\begin{equation*}0\to\mathcal I_h(K_X+m_0L)\to K_X+m_0L\to ({\mathcal O_X}/{\mathcal I_h}) (K_X+m_0L)\to 0.\end{equation*}
The relevant part of the pertaining long exact sequence reads:
\begin{equation*}\Gamma (X,K_X+m_0L)\to\Gamma(V(\mathcal I_h),({\mathcal O_X}/{\mathcal I_h})(K_X+m_0L)) \to 0,\end{equation*}
which implies by virtue of the second property that
\begin{equation*}\Gamma (X,K_X+m_0L)\stackrel{\text{restr.}}{\longrightarrow}\Gamma(\{x\},{\mathcal O}_{\{x\}}(K_X+m_0L)) \to 0,\end{equation*}
meaning that the restriction map to $\{x\}$ is surjective, which is what we intended to prove. Note that, since $L$ is ample, there trivially exists a metric with the two aforementioned properties for every line bundle $mL$ with $m\geq m_0$ (just multiply the metric for $m_0L$ by the $(m-m_0)$-th power of a smooth positive metric of $L$).\par
In both \cite{AS} and \cite{H}, the sought-after metric $h$ is produced by an inductive method. First, here is the key statement proved in sections 7--9 of \cite{AS}. The cornerstone of its proof is a clever application of the theorem of Ohsawa and Takegoshi on the extension of $L^2$ holomorphic functions \cite{OT}. Note that, in contrast to \cite{AS}, we are only concerned with freeness and not point separation, so we can do without the complicated formulations found there.
\begin{proposition}[\cite{AS}]\label{asind}
Let $d$ be an integer with $1\leq d\leq n-1$. Let $k_d$ be a positive rational number, and let $h_d$ be a singular metric of the line bundle $k_dL$. Assume that $x\in  V(\mathcal I_{h_d})$ and $x \not \in V(\mathcal I_{(h_d)^\gamma})$ for $\gamma <1$. Moreover, assume that the dimensions of those components of $V(\mathcal I_{h_d})$ which contain $x$ do not exceed $d$. Then there exist integers $d',k_{d'}$ with $0\leq d'<d$ and $k_d<k_{d'}<k_d+d+\varepsilon$ ($\varepsilon$ denotes a positive rational number which can be chosen to be arbitrarily small) and a singular metric $h_{d'}$ of $k_{d'}L$ such that $h_{d'}$ possesses the same properties as $h_d$, but with $d$ and $k_d$ replaced by $d'$ and $k_{d'}$.
\end{proposition}
Secondly, the key statement of \cite{H} is the following proposition. It is stated in such a way that it unites \cite{H}, Proposition 3.2 (the inductive statement), and \cite{H}, Corollary 4.6 (the multiplicity bound), into one ready-to-use statement. In its proof, the use of the aforementioned $L^2$ extension theorem is avoided by an explicit bound on the multiplicity of the minimal centers occurring in the inductive procedure.
\begin{proposition}[\cite{H}]\label{helmkeind}
Let $d$ be an integer with $1\leq d\leq n-1$. Let \begin{equation*} L^n> n^n\end{equation*}
and 
\begin{equation*}L^{\tilde d}.Z\geq m_x(Z)\cdot n^{\tilde d}\end{equation*}
for all subvarieties $Z \subset X$ such that $x\in Z,\ d\leq \tilde d= \dim Z \leq n-1$ and multiplicity $m_x(Z)\leq{{n-1}\choose{\tilde d-1}}$ at $x$. Then there exists an integer $0 \leq d'< d$, a rational number $0<c<1$ and an effective $\Q$-divisor $D$ such that $D$ is $\Q$-linearly equivalent to $cL$, the pair $(X,D)$ is log canonical at $x$ and the minimal center of $(X,D)$ at $x$ is of dimension $d'$.
\end{proposition}
Let us briefly recall the definitions of some of the terms occurring in Proposition \ref{helmkeind}. First of all, for a pair $(X,D)$ of a variety $X$ and a $\Q$-divisor $D$, an {\it embedded resolution} is a proper birational morphism $\pi: Y\to X$ from a smooth variety $Y$ such that the union of the support of the strict transform of $D$ and the exceptional divisor of $\pi$ is a normal crossing divisor. With this basic definition, we continue.
\begin{definition}
Let $X$ be a normal variety and $D=\sum_id_iD_i$ an effective $\Q$-divisor such that $K_X+D$ is $\Q$-Cartier. If $\pi: Y\to X$ is a birational morphism (in particular, an embedded resolution of the pair $(X,D)$), we define the {\it discrepancy divisor} of $(X,D)$ under $\pi$ to be
\begin{equation*}\sum_jb_jF_j:=K_Y-\pi^*(K_X+D).\end{equation*}
The pair $(X,D)$ is called {\it log canonical} (resp.~{\it Kawamata log terminal}) at $x$, if there exists an embedded resolution $\pi$ such that $b_j\geq -1$ (resp.~$b_j>-1$) for all $j$ with $x \in \pi(F_j)$. Moreover, a subvariety $Z$ of $X$ containing $x$ is said to be a {\it center of a log canonical singularity} at $x$, if there exists a birational morphism $\pi: Y\to X$ and a component $F_j$ with $\pi(F_j)=Z$ and $b_j\leq -1$.
\end{definition}\par
It follows from Shokurov's connectedness lemma in \cite{Sh} that the intersection of two centers of a log canonical singularity is again a center of a log canonical singularity (for a proof, see \cite{K}). Thus there exists a unique minimal center of a log canonical singularity at $x$ with respect to the inclusion of subvarieties on $X$.\par
In order to make use of Proposition \ref{helmkeind} for our purposes, we derive from its conclusion a statement about the existence of a certain singular metric.
\begin{proposition}\label{translate}
Let $(X,D)$ be a pair of a smooth projective variety $X$ and an effective $\Q$-divisor $D$. Let $x\in X$ be an arbitrary but fixed point. Assume that the pair $(X,D)$ is a log canonical at $x$ with its minimal center at $x$ being non-empty. Let $0<c<1$ be a rational number such that $D$ is $\Q$-linearly equivalent to $cL$. Then there exists a singular metric $h_D$ and a rational number $c'$ (which can be chosen to be arbitrarily close to c) such that $h_D$ is a metric of $c'L$, $x\in V(\mathcal I_{h_D})$ and $V(\mathcal I_{h_D})$ is contained in the minimal center of $(X,D)$ at $x$ in a neighborhood of $x$. Moreover, $x\not\in V(\mathcal I_{(h_D)^\gamma})$ for $\gamma < 1$.
\end{proposition}
\begin{proof}
Let $s$ be a multivalued holomorphic section of $cL$ whose $\Q$-divisor is $D$. This means that for some positive integer $p$ with $cp$ being an integer, the $p$-th power of $s$ is the canonical holomorphic section of $pcL$ with divisor $pD$. Let $Z$ denote the minimal center of $(X,D)$ and $\pi:Y\to X$ a log resolution of $(X,D)$ with discrepancy divisor $\sum_jb_jF_j$. We choose $\pi$ such that there exists at least one index $j_0$ with $b_{j_0}=-1$ and $\pi(F_{j_0})=Z$. Furthermore, we set $\sum_j\delta_jF_j:=\pi^*(D)$.\par
Since $L$ is ample, we can choose a finite number of multivalued holomorphic sections $s_1,\ldots,s_q$ of $L$ whose common zero set is exactly $Z$. Let $\delta_{i,j}$ denote the vanishing order of $\pi^*s_i$ along $F_j$ at a generic point of $F_j$. If we set $\delta:=\min\{\delta_{i,j_0}:i=1,\ldots,q\}$, then $\delta>0$ holds because all $s_i$ vanish on $Z$.\par
For small positive rational numbers $\varepsilon,\varepsilon'<1$, we define the following singular metric of, say, $\tilde cL$:
\begin{equation*}\tilde h_D:=\frac 1 {|s|^{2(1-\varepsilon)}}\frac 1 {(\sum_{i=1}^q|s_i|^2)^{\varepsilon'}}.\end{equation*}
Whatever the choice of $\varepsilon,\varepsilon'$ may be, $\tilde h_D$ is locally integrable outside of $Z$ in a small neighborhood of $x$. Here is how to choose $\varepsilon,\varepsilon'$ in order to make $\tilde h_D$ not integrable at $x$. In a small neighborhood $U$ of $x$, the integrability of $\tilde h_D$ is equivalent to the integrability of
\begin{equation*}\pi^*\tilde h_D|\text{Jac}(\pi)|^2=\frac 1 {|s\circ\pi|^{2(1-\varepsilon)}}\frac 1 {(\sum_{i=1}^q|s_i\circ\pi|^2)^{\varepsilon'}}|\text{Jac}(\pi)|^2\end{equation*}
over every small open subset $W$ of $\pi^{-1}(U)$. Note that $|\text{Jac}(\pi)|^2$ can only be defined locally, and over $W$ we take it to be the quotient 
\begin{equation*}\frac {\pi^*(\omega_U\wedge\bar{\omega}_U)}{\omega_W\wedge\bar{\omega}_W},\end{equation*}
where $\omega_U,\omega_W$ are arbitrary but fixed nowhere vanishing local holomorphic $n$-forms on $U$ and $W$, respectively. As we continue, we observe that there exists a small open subset $W$ of $\pi^{-1}(U)$ such that $W\cap F_{j_0}\neq\emptyset$ and $\pi^*\tilde h_D|\text{Jac}(\pi)|^2$ has a pole along $W\cap F_{j_0}$, with its order at a generic point of $W\cap F_{j_0}$ being
\begin{equation*}b_{j_0}+\varepsilon\delta_{j_0}-\varepsilon'\delta.\end{equation*}
This number equals $-1$ if we choose $\varepsilon$ arbitrarily and set $\varepsilon':=\frac 1 {\delta}\varepsilon\delta_{j_0}$. We conclude that, with these choices for $\varepsilon$ and $\varepsilon'$, $\tilde h_D$ is not integrable at $x$.\par
Finally, we set $h_D:=(\tilde h_D)^r$ with $r:=\min\{\rho:0<\rho\leq 1,\ (\tilde h_D)^\rho \text{ is not integrable at } x\}$ to obtain the desired singular metric for $c'L$. Notice that if we let $\varepsilon\to 0$, then $r\to 1$, $\tilde c \to c$ and $c'\to c$.
\end{proof}
Now we are in a position to prove our Main Theorem.
\begin{proof} [Proof of the Main Theorem] 
Fix $x\in X$. Our goal is to prove that, if $m_0$ is the smallest integer no less than
\begin{equation*}(e+\frac 1 2)n^\frac 4 3 + \frac 1 2 n^\frac 2 3 +1,\end{equation*}
there exists a singular metric $h$ of the line bundle $m_0L$ such that the two properties listed at the beginning of this section are satisfied. As was explained before, this is all that is necessary to prove the Main Theorem.\par
Let $a$ be the smallest integer which is no less than $e\,n^{\frac 4 3}$. Let $d_0$ be the integral part of $n^{\frac 2 3}$. According to Lemma \ref{coeffestim}, we have
\begin{equation*}
{n-1\choose \tilde d-1}^{\frac 1 {\tilde d}}n\leq e\,\frac{n^2}{\tilde d}.
\end{equation*}
for all integers $\tilde d$ with $1\leq \tilde d \leq n-1$.
Furthermore,
\begin{equation*}
e\,\frac{n^2}{\tilde d}\leq e\,n^{\frac 4 3}\leq a
\end{equation*}
for $\tilde d\geq n^{\frac 2 3}$. Thus we can use Proposition \ref{helmkeind} to produce an effective $\Q$-divisor $D$ such that $D$ is $\Q$-linearly equivalent to $caL$ for some $0<c<1$, the pair $(X,D)$ is log canonical at $x$ and its minimal center at $x$ is of dimension $d'$ for some integer $d'$ with $0\leq d' \leq d_0$.
By Proposition \ref{translate}, this translates into the existence of a singular metric $h_1$ of $c'aL$ $(0<c'<1)$ such that $x\in V(\mathcal I_{h_1})$, $x\not \in V(\mathcal I_{(h_1)^\gamma})$ for $\gamma < 1$ and the dimensions of those components of $V(\mathcal I_{h_1})$ which contain $x$ do not exceed $d'$.\par
From this point onwards, we can use the method of \cite{AS} in the form of Proposition \ref{asind} to produce inductively a singular metric $h_2$ such that $V(\mathcal I_{h_2})$ is isolated at $x$. If the constructed metric $h_2$ is a metric for, say, $kL$, then  
\begin{equation*}k\leq c'a+1+2+\ldots+d_0+\varepsilon_1+\varepsilon_2+\ldots+\varepsilon_{d_0}.\end{equation*}
Since the $\varepsilon_i$ can be chosen to be arbitrarily small positive rational numbers and since $c'<1$, we can assume that
\begin{equation*}c'a+1+2+\ldots+d_0+\varepsilon_1+\varepsilon_2+\ldots+\varepsilon_{d_0}<a+1+2+\ldots+d_0.\end{equation*}\par
In order to obtain a metric of $m_0L$ with the additional property that its curvature current dominates a positive definite smooth $(1,1)$-form on $X$ in the sense of currents, we can simply multiply $h_2$ by the $(m_0-k)$-th power of a smooth positive metric of $L$ to obtain the desired metric $h$ of $m_0L$. Note that $m_0-k$ is a positive number because
\begin{eqnarray*}
m_0-k&>&m_0-(a+1+2+\ldots+d_0)\\
&\geq& m_0-(e\,n^{\frac 4 3}+1+1+2+\ldots+d_0)\\
&=& m_0-(e\,n^{\frac 4 3}+1+\frac 1 2 d_0(d_0+1))\\
&\geq&m_0-(e\,n^{\frac 4 3}+1+\frac 1 2 n^{\frac 2 3}(n^{\frac 2 3}+1))\\
&=&m_0-((e+\frac 1 2)n^\frac 4 3 + \frac 1 2 n^\frac 2 3 +1)\geq0.
\end{eqnarray*}
The proof of the Main Theorem is now complete.
\end{proof}

\section{Applications}
As was indicated before, not much is known about the very ampleness part of the Fujita conjecture. The theorems and techniques mentioned in the previous sections do not seem to be directly applicable to it. However, Angehrn and Siu \cite{AS} were able to prove the following weaker analog to the very ampleness part of Fujita's conjecture, in which they assume that $L$, in addition to being ample, is also base point free. Their result improves on previous results of Ein, K\"uchle and Lazarsfeld \cite{EKL} and Kollar \cite{Ko}.
\begin{theorem}[\cite{AS}]\label{freetoample}
Let $L$ be an ample line bundle over a compact complex manifold $X$ of complex dimension $n$ such that $L$ is free. Let $A$ be an ample line bundle. Then $(n+1)L+A+K_X$ is very ample.
\end{theorem}
In conjunction with our Main Theorem, Theorem \ref{freetoample} can be readily applied to the case of an ample canonical line bundle in order to give the following effective statement on pluricanonical embeddings. As far as the author knows, this is best effective statement on pluricanonical embeddings currently on hand. Note that Fujita's conjecture indicates that the statement of the corollary should hold true for any integer $m\geq n+3$.
\begin{corollary}\label{effpluri}
If $X$ is a compact complex manifold of complex dimension $n$ whose canonical bundle $K_X$ is ample, then $mK_X$ is very ample for any integer $m\geq (e+\frac 1 2)n^\frac 7 3+\frac 1 2 n^\frac 5 3 + (e+\frac 1 2)n^\frac 4 3 + 3n+ \frac 1 2 n^\frac 2 3+5$.
\end{corollary}
\begin{proof}
Let $m_0$ be the smallest integer no less than $(e+\frac 1 2)n^\frac 4 3 + \frac 1 2 n^\frac 2 3 +1$. According to our Main Theorem, $m_0K_X+K_X=(m_0+1)K_X$ is base point free (and, of course, ample). Thus we can apply Theorem \ref{freetoample} with $L=(m_0+1)K_X$ and $A=K_X$ to obtain that $mK_X$ is very ample for any integer $m\geq (n+1)(m_0+1)+2$. A simple estimate yields the following upper bound for $(n+1)(m_0+1)+2$:
\begin{eqnarray*}
&&(n+1)(m_0+1)+2\\
&=&m_0(n+1)+n+3\\
&\leq& ((e+\frac 1 2)n^\frac 4 3 + \frac 1 2 n^\frac 2 3 +1+1)(n+1)+n+3\\
&=&(e+\frac 1 2)n^\frac 7 3+\frac 1 2 n^\frac 5 3 + (e+\frac 1 2)n^\frac 4 3 + 3n+ \frac 1 2 n^\frac 2 3+5, \text{ q.e.d.}
\end{eqnarray*} 
\end{proof}
Finally, we remark that our effective statement on pluricanonical embeddings can be used to sharpen the best known bound for the number of dominant holomorphic maps from a fixed compact complex manifold with ample canonical bundle to any variable compact complex manifold with big and numerically effective canonical bundle.
\begin{acknowledgments}
The author is supported by a doctoral student fellowship of the Studienstiftung des deutschen Volkes (German national merit foundation). The results published in this article were obtained as part of the author's research for his Ph.D. dissertation at the Ruhr-Universit\"at Bochum under the auspices of Professor Alan T. Huckleberry. It is a great pleasure to thank Professor Yum-Tong Siu for introducing the author to the field of effective algebraic geometry and for numerous discussions on the subject, which took place during visits to the Mathematics Department of Harvard University and the Institute of Mathematical Research at the University of Hong Kong. Special thanks go to Professor Ngaiming Mok for making an extended visit to Hong Kong possible. The author is a member of the Forschungs\-schwerpunkt ``Globale Methoden in der komplexen Geometrie'' of the Deutsche Forschungsgemeinschaft.
\end{acknowledgments}

\end{document}